\documentclass[11pt,reqno]{amsart}

\usepackage{amsmath}
\usepackage{amsthm}
\usepackage{amsfonts}
\usepackage{amssymb}
\usepackage{color}
\usepackage{graphicx}
\usepackage[hidelinks]{hyperref}
\usepackage[margin=1.0in]{geometry}
\usepackage{todonotes}
\usepackage{verbatim}
\usepackage{comment}
\numberwithin{equation}{section}
\usepackage{enumitem}
\setlist{nosep}
\usepackage{caption}
\usepackage{subcaption}
\usepackage{epstopdf}
\usepackage[numbers]{natbib}
\usepackage{longtable}

\definecolor{skyblue}{rgb}{0.85,0.85,1}

\newtheorem{rem}{Remark}

\newtheorem{define}{Definition}
\newtheorem{thm}{Theorem}
\numberwithin{lemma}{section}
\numberwithin{prop}{subsection}
\numberwithin{theorem}{section}
\numberwithin{cor}{section}
\numberwithin{conj}{section}
\numberwithin{rem}{section}

\newcommand{\bbC}{\mathbb{C}}

\newcommand{\bbR}{\mathbb{R}}

\newcommand{\bbZ}{\mathbb{Z}}
\newcommand{\vp}{\varphi}

\newcommand{\p}{\partial}

\DeclareMathOperator{\Gr}{Gr}

\newcommand{\eps}{\epsilon}
\newcommand{\Maslov}{\mathrm{Maslov}(\vp_\eps)}

\begin{document}

\title[Bifurcation to instability through the lens of the Maslov index]{Bifurcation to instability through the lens of the Maslov index}
\author[P. Cornwell]{Paul Cornwell}
\address{Johns Hopkins University Applied Physics Lab, 11100 Johns Hopkins Road, Laurel, MD 20723-6099}
\email{paul.cornwell@jhuapl.edu}

\address{Department of Mathematics, UNC Chapel Hill, Phillips Hall CB \#3250, Chapel Hill, NC 27516}
\author[C.\,K.\,R.\,T.\ Jones]{Christopher K.\,R.\,T.\ Jones}
\email{ckrtj@email.unc.edu}
\author[C. Kiers]{Claire Kiers}
\email{cekiers@live.unc.edu   }

\begin{abstract}
The Maslov index is a powerful tool for assessing the stability of solitary waves. Although it is difficult to calculate in general, a framework for doing so was recently established for singularly perturbed systems \cite{CJ18}. In this paper, we apply this framework to standing wave solutions of a three-component activator-inhibitor model. These standing waves are known to become unstable as parameters vary. Our goal is to see how this established stability criterion manifests itself in the Maslov index calculation. In so doing, we obtain new insight into the mechanism for instability. We further suggest how this mechanism might be used to reveal new instabilities in singularly perturbed models.
\end{abstract}

\maketitle

\tableofcontents

\section{Introduction}

Many physical systems conserve energy. In such systems, the notion that states which minimize energy are stable is a fairly reliable heuristic. Solitary waves, for example, often arise as critical points of an energy functional. One checks whether a wave is a minimizer by analyzing the second derivative of the energy, which happens to be the linearization about the wave. The existence of unstable spectrum therefore indicates that there are directions in which the energy decreases. For  variational problems in general, the \emph{Morse index} is defined to be the dimension of the maximal subspace on which the Hessian of the energy evaluated at a critical point is negative definite. 

Computing the Morse index amounts to an infinite-dimensional eigenvalue problem, which is typically difficult to solve. The celebrated Morse Index Theorem \cite[\S 15]{Milnor} is a valuable tool in the special case where the functional is the energy of a path, and critical points are geodesics on a manifold $M$. It states that the Morse index, evaluated at a critical path $\gamma(t)$, is equal to the number of conjugate points along $\gamma$. Without diving into definitions, the crux of this result is that the Morse index is determined by how $\gamma$ is situated in the tangent bundle $TM$; there is no need to analyze the ``spectrum" of the Hessian explicitly.

The subject of this paper is an adaptation--in fact, a generalization--of the Morse Index Theorem to the context of stability of solitary waves. The number of conjugate points for a critical path is replaced by the Maslov index of the wave, which is an intersection number assigned to curves of Lagrangian planes. One can show that the Maslov index counts (or at least gives a lower bound for) real, unstable eigenvalues. Moreover, it is computed by fixing the spectral parameter $\lambda=0$, which yields the equation of variations for the wave. This important fact is the key to extracting stability information from phase space geometry, which we discuss later.

The equality of the Morse and Maslov indices has been worked out in a number of settings in recent years, e.g. \cite{CH,CH14,corn,HLS16,JLS17,JLM13}. The bigger challenge is arguably calculating the Maslov index, which has proven difficult to do. In the spirit of the Morse Index Theorem, two schools have emerged with techniques for doing so. The first is based on the calculus of variations, owing to Chen and Hu \cite{CH, CH14}. They use the Maslov index to show that the Morse index is $0$ for energy minimizers. The ``dirty work" of the calculation is then to construct an admissible class of functions and find a minimizer.  This strategy was employed to prove the existence and (in)stability of standing waves in a doubly-diffusive FitzHugh-Nagumo system \cite{CC12, CH14}.

The second school strives to calculate or estimate the Maslov index directly by locating intersections. This approach has produced instability results for generic standing waves of gradient reaction-diffusion systems \cite{BCJLMS17}, as well as (in)stability results for various standing and traveling waves \cite{BJ,CDB09,CDB_multi,CJ18,Jo88}. The main technical tool for this approach is the crossing form of Robbin and Salamon \cite{RS93}. A conjugate point corresponds to the intersection of a curve of Lagrangian planes with a codimension one set (the ``singular cycle"), and the crossing form determines the contribution to the Maslov index at a conjugate point.

We now briefly describe the challenge of calculating the Maslov index directly. Let $L$ be the operator obtained by linearizing about a solitary wave. The eigenvalue equation $Lp=\lambda p$ can be cast as a non-autonomous dynamical system on $\bbR^{2n}$ (for $\lambda\in\bbR$), where $n$ is the number of components. The curve of interest in the Maslov index calculation is the $n$-dimensional subspace of solutions to $Lp=0$ which decay at $-\infty$, called the unstable bundle. Due to translation invariance, the derivative of the wave is everywhere part of this space. In the scalar case, this is the only solution, and conjugate points can be related to zeros of the velocity (from which Sturm-Liouville theory follows almost immediately). It is the presence of other solutions in the higher-dimensional case that makes the calculation difficult. Finding these solutions is tantamount to solving a linear, non-autonomous equation on the real line.

In \cite{CJ18}, two of the authors of this paper established a framework for calculating the Maslov index in singularly perturbed equations using geometric singular perturbation theory (GPST). The strategy is to use the fact that the unstable \emph{bundle} is everywhere tangent to the unstable \emph{manifold} containing the wave in phase space. Using Fenichel theory \cite{Fen79} and subsequent developments such as the Exchange Lemma \cite{JK94}, one can accurately determine the orientation of the unstable manifold as it evolves along the wave.\footnote{The reader must decide whether tracking the unstable manifold should be done with \cite{JK94} or without \cite{Bru96} the use of differential forms.} Thus it is possible to calculate the Maslov index without having to solve the linear system explicitly. As a proof of concept, this framework was applied to show that fast traveling waves for a FitzHugh--Nagumo equation (with equal diffusion rates) are stable.

The aim of this paper is to apply the framework of \cite{CJ18} to a system of three reaction-diffusion equations (\ref{PDE model slow}) introduced by Schenk et al \cite{Schenk}. Using GSPT, Kaper et al showed that (\ref{PDE model slow}) supports a multitude of interesting standing and traveling wave solutions \cite{DvHK_ex}. They subsequently obtained stability results using a fast-slow Evans function decomposition \cite{vHDK_stab}. These results were then reproved and expanded upon by van Heijster et al using a variational approach \cite{vHCNT}. The latter stability proof uses the Maslov index in the manner of the first school described above. We will use a conjugate point-based Maslov index calculation to re-derive the known stability result in the case of standing single pulses.

We have several objectives in this work. The first is to demonstrate the robustness of the calculation method developed in \cite{CJ18}. In particular, we prove its utility in systems of dimension greater than two--the simplest non-trivial case. Second, the pulses that we study can be stable or unstable depending on the model parameters. By tuning the parameters, we can therefore see exactly how the instability manifests itself in the Maslov index calculation. This stands in contrast to the calculation in \cite{CJ18}, where the waves are stable regardless of the parameters (at least locally). Using this new example, we will argue that the Maslov index can be used to identify and manipulate mechanisms for (in)stability in singularly perturbed systems.

Finally, we believe that the calculation carried out herein will be of interest to GSPT itself. Since the original work of Fenichel \cite{Fen79}, the two most important advances to the geometric theory of singularly perturbed systems are the Exchange Lemma and geometric desingularization--the ``blowup method" \cite{Dum93}. The former is useful when considering passage near a normally hyperbolic critical manifold, whereas the latter is used to study dynamics at a point where normal hyperbolicity is lost. We shall see later that the Maslov index calculation requires a hybrid approach. Although the system we analyze has no fold points, the fast-slow transitions in the tangent bundle above the wave are critical. As such, it is necessary to zoom in on the exact point in phase space where this transition occurs. Instead of studying the dynamics near these points on a sphere--as one would do in the blowup method--we aim to understand the dynamics of the equation of variations on the Lagrangian Grassmannian. An interesting extension of this work would be to calculate the Maslov index of a wave which passes through a fold point, e.g. \cite{CRS}.

The rest of this paper is organized as follows.  In section 2, we state the model equations and emphasize the fast-slow structure. We then describe the standing pulse solutions of interest by breaking them into fast and slow components. In section 3, we lay out the eigenvalue problem and define the Maslov index of the pulses. In section 4, we carefully compute the Maslov index and show how an instability can appear when the model parameters are changed. Finally, we conclude in section 5 by explaining the insight gained from this calculation as well as how it might apply to other singularly perturbed systems.

\section{A 3-component activator-inhibitor model}

In this work, we consider the following three-component system of reaction-diffusion equations:
\begin{equation}\label{PDE model slow}
\begin{aligned}
U_t & = \eps^2U_{xx}+U-U^3-\eps(\alpha V+\beta W+\gamma)\\
\tau V_t & = V_{xx}+U-V\\
\theta W_t & =D^2W_{xx}+U-W.
\end{aligned}
\end{equation} We assume that $\tau,\theta>0$, $D>1$, $0<\eps\ll 1$, and $\alpha,\beta,\gamma\in\bbR$. We further assume that all parameters are $O(1)$ in $\eps$. This model has its roots in the high-level study of gas-discharge dynamics. For more background on the physical aspects of the problem, we refer the reader to \cite{vHCNT}. In \cite{DvHK_ex,vHDK_stab}, geometric singular perturbation theory was used to prove the existence and stability of various standing and traveling waves for (\ref{PDE model slow}). These results were then revisited and expanded in \cite{vHCNT} by combining the GSPT analysis with an action functional approach.

We will consider standing pulse solutions of (\ref{PDE model slow}), which are time-independent, localized structures. By setting $U_t=V_t=W_t=0$ and introducing the variables $P=\eps U_x, Q=V_x,$ and $R=DW_x$, such solutions are seen to be homoclinic orbits for the standing-wave ODE
\begin{equation}
\begin{split}\label{standing wave ODE slow}
\epsilon U' & = P \\
\epsilon P' & = -U+U^3 + \epsilon(\alpha V + \beta W + \gamma) \\
V' & = Q \\
Q' & = V-U \\
W' & = \frac{1}{D} R \\
R' & = \frac{1}{D} (W-U)
\end{split},
\end{equation} where $' = \frac{d}{dx}.$ (\ref{standing wave ODE slow}) is formulated on the ``slow" timescale. We will also use the ``fast" version of (\ref{standing wave ODE slow}), obtained by setting $\eps\xi=x.$ Denoting $\dot \ = \frac{d}{d\xi}$, we rewrite (\ref{standing wave ODE slow}) as
\begin{equation}
\begin{split}\label{standing wave ODE fast}
\dot U & = P \\
\dot P & = -U + U^3 + \epsilon(\alpha V + \beta W + \gamma) \\
\dot V & = \epsilon Q \\
\dot Q & = \epsilon (V-U) \\
\dot W & = \frac{\epsilon}{D} R \\
\dot R & = \frac{\epsilon}{D} (W-u)
\end{split}.
\end{equation} When $\eps>0$, (\ref{standing wave ODE slow}) and (\ref{standing wave ODE fast}) define the same dynamics. However, the limiting systems obtained by sending $\eps\rightarrow 0$ are distinct. The fast subsystem (\ref{standing wave ODE fast}) is two-dimensional with $V,Q,W,$ and $R$ acting as parameters. On the other hand, (\ref{standing wave ODE slow}) is a differential-algebraic equation where the dynamics are four-dimensional and restricted to the set \begin{equation}\label{slow mfd defn}
    M_0 = \{(U,P,V,Q,W,R)\in \bbR^6: P=0, U\in \{-1,0,1\}\}.
\end{equation} $M_0$ is called the \emph{critical manifold}. Observe that $M_0$ is precisely the set of critical points for (\ref{standing wave ODE fast}) with $\eps=0$. 

The goal of GSPT is to construct solutions to (\ref{standing wave ODE fast}) for small $\eps>0$ by gluing together segments from the fast and slow flows. The technical challenge is to prove that these singular orbits persist when $\eps>0$. Since this theory is well understood (at least as applied to the system at hand), we will simply describe the $\eps=0$ object, which is all that is needed to compute the Maslov index. The interested reader can find more detail on the existence proofs in \cite[\S 2.1-2.3]{DvHK_ex}, or on GSPT more generally in \cite{Jones_GSP,Kuehn15}.

\subsection{Standing single pulses}

As mentioned earlier, there are many permanent structures hiding in this model.  Since the objective of this paper is to observe the bifurcation to instability through the lens of the Maslov index, we will focus on the simplest structure that exhibits this behavior--the standing single pulse.

In phase space, a solitary pulse is a homoclinic orbit to a fixed point. For $\eps>0$, once sees that (\ref{standing wave ODE fast}) has three fixed points, which are distinguished by the $U$ component. The wave we consider is homoclinic to the smallest value of $U$, which we see from (\ref{standing wave ODE fast}) and (\ref{slow mfd defn}) is \begin{equation}\label{u_smallest} U_{\eps}^-:=-1+O(\eps).\end{equation}
The other components of the fixed point follow directly from (\ref{u_smallest}). We define the fixed point \begin{equation}
    X_\eps^-:=(U_{\eps}^-,0,U_{\eps}^-,0,U_{\eps}^-,0)
\end{equation} of (\ref{standing wave ODE fast}) as the rest state of the wave.
The homoclinic orbit itself will be $O(\eps)$ close to a singular version consisting of five segments: three slow and two fast. We denote the singular orbit $\vp_0(x)\in\bbR^6$, where the subscript refers to $\eps=0$. To describe $\vp_0$ in more detail, we will need the two-dimensional fast system \begin{equation}
\begin{split}\label{fast reduced system}
\dot U & = P \\
\dot P & = -U + U^3,
\end{split}
\end{equation} and the four-dimensional slow reduced system
\begin{equation}
\begin{split}\label{slow reduced system}
V' & = Q \\
Q' & = V-U \\
W' & = \frac{1}{D} R \\
R' & = \frac{1}{D} (W-U).
\end{split}
\end{equation}

The two fast orbits are heteroclinic connections from $(U,P)=(-1,0)$ to $(U,P)=(+1,0)$ and back again. Observe that (\ref{fast reduced system}) is Hamiltonian with $$H(U,P)=P^2/2+U^2/2-U^4/4,$$ so we can solve $H(U,P)=H(-1,0)$ to see that the heteroclinic connections are given by \begin{equation}\label{up_hetero}
P=\pm\sqrt{\frac{1}{2}-U^2+\frac{U^4}{2}}=\pm\frac{1}{\sqrt{2}}(1-U^2).    
\end{equation} (The positive branch goes from $(-1,0)$ to $(1,0)$.)

For the slow segments, $U$ is now a parameter. Setting $U=-1$ in (\ref{slow reduced system}), we see that the dynamics are linear, except with the origin shifted to \begin{equation}
    (V, Q, W, R) = (-1,0,-1,0):=X_-.
\end{equation} Moreover, the $(V,Q)$ and $(W,R)$ equations decouple, so it is easy to compute the (2D) stable and unstable manifolds of $X_-$: \begin{equation}
\begin{split}
    W^s(X_-) & = \{(V,Q,W,R):Q=-(V+1), R = -(W+1)\} \\
    W^u(X_-) & = \{(V,Q,W,R):Q=V+1, R = W+1\}.
\end{split}
\end{equation} 

Likewise, for $U=+1$, (\ref{slow reduced system}) is linear with the origin shifted to \begin{equation}
    (V, Q, W, R) = (1,0,1,0):=X_+.
\end{equation} The stable and unstable manifolds of $X_+$ are \begin{equation}
\begin{split}
    W^s(X_+) & = \{(V,Q,W,R):Q=-(V-1), R = -(W-1)\} \\
    W^u(X_+) & = \{(V,Q,W,R):Q=V-1, R = W-1\}.
\end{split}
\end{equation}

We can now describe how the segments fit together to form $\vp_0$. First, since these manifolds will play a prominent role later, we officially define the relevant subsets of the critical manifold:
\begin{equation}
\begin{split}
    M_0^- & = \{(U,P,V,Q,W,R):U=-1,P=0\} \\
    M_0^+ & = \{(U, P, V, Q, W, R):U=1,P=0\}.
\end{split}
\end{equation} To simplify notation, we also append the appropriate $(U,P)$ coordinates to $W^{u/s}(X_\pm)$ so that $W^{u/s}(X_-)\subset M_0^-$, and $W^{u/s}(X_+)\subset M_0^+$. The two fast segments, as mentioned above, are heteroclinic connections between $M_0^-$ and $M_0^+$. Two of the slow segments will turn out to be trajectories inside of $W^{u/s}(X_-)$; the first segment leaves $X_-$ along $W^u(X_-)$, and the fifth segment returns to $X_-$ along $W^s(X_-)$. 

\begin{figure}
\centering
\includegraphics[scale = 0.5]{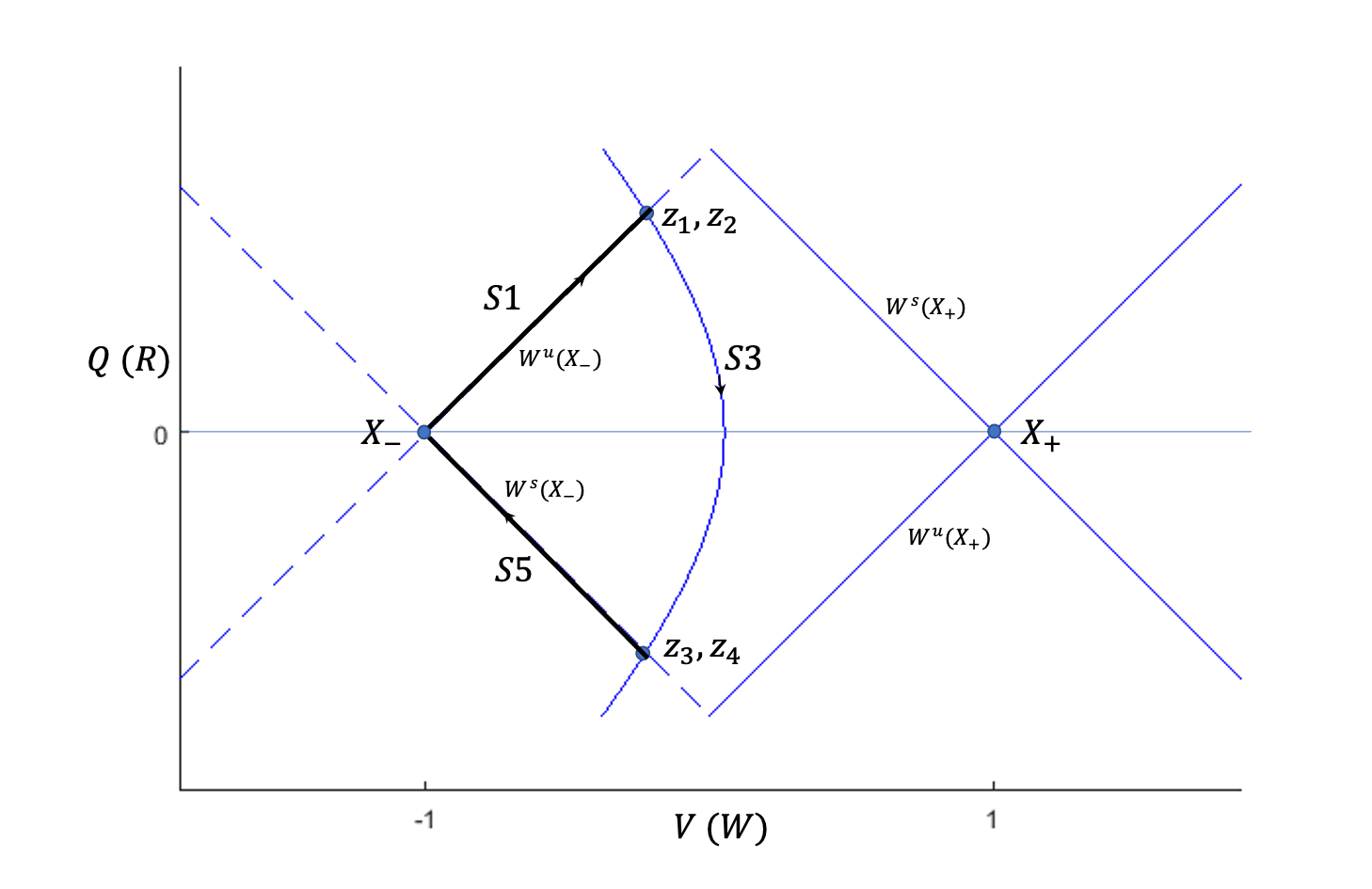}
\caption{Projection of slow flows on $M_0^+$ and $M_0^-$ onto $VQ-$space (or $WR-$space). Segment 3 (on $M_0^+$) serves to carry the slow coordinates from $W^u(X_-)$ to $W^s(X_-)$.}
\label{slow_flow_projection}
\end{figure}

The third segment is the part of the slow flow on $M_0^+$. Since the slow coordinates do not change along the fast jumps, the role of this segment is to flow the slow coordinates from $W^u(X_-)$ to $W^s(X_-)$. The only way for this to happen is for the $V$ and $W$ coordinates of the jump-off point to be negative. This is equivalent to the jump-off occurring prior to the intersection of $W^u(X_-)$ and $W^s(X_+)$ projected onto the slow coordinates. (See Figure \ref{slow_flow_projection}.) Since (\ref{slow reduced system}) is linear, segment 3 lies on the product of two hyperbolas--one each in $VQ-$ and $WR-$space.

At this point, we have completely described the singular orbit qualitatively. However, one detail remains: the values of $V$ and $W$ at the jump-off points. We write $z_i,i\in\{1,2,3,4\},$ for the point in $\bbR^6$ where segment $i$ meets segment $i+1$. From Figure \ref{slow_flow_projection}, it is clear that specifying $V$ and $W$ at $z_1$ determines all four variables at each $z_i$. It turns out that fixing these values is an $\eps>0$ consideration based on the change in the Hamiltonian along the fast jump. We refer the reader to \cite{DvHK_ex} for the derivation and instead just state the existence result for the standing pulses.

\begin{thm}[Theorem 2.1 of \cite{DvHK_ex}, Theorem 1 of \cite{vHCNT}]\label{thm_existence}
Let $(\alpha, \beta, \gamma, D)$ be such that \begin{equation}\label{jump_cond_x}
    \alpha e^{-2x}+\beta e^{-2x/D}=\gamma
\end{equation} has a positive solution $x=x^*$. Then, for $\eps>0$ sufficiently small, (\ref{standing wave ODE fast}) possess a homoclinic orbit $\vp_\eps(x)$, whose $U-$coordinate approaches $U_\eps^-$ defined by (\ref{u_smallest}) as $x\rightarrow\pm\infty$. Moreover, the image of $\vp_\eps(x)$ is $O(\eps)$-close to $\vp_0$ with jump-off point \begin{equation}\label{jump_off_coords}
    z_1=(-1,0,-e^{-2x^*},1-e^{-2x^*},-e^{-2x^*/D},1-e^{-2x^*/D}).
\end{equation}
\end{thm}

\begin{rem}
    Depending on the parameters, equation (\ref{jump_cond_x}) can have zero, one, or two solutions. If there are zero solutions, no single pulse waves exist.  If there are two solutions, then two waves exist with different jump-off values.
\end{rem}

\begin{rem}\label{rem slow mfd jump cond}
    Comparing (\ref{jump_cond_x}) and (\ref{jump_off_coords}), we see that the jump-off condition can be rewritten in terms of the slow coordinates as \begin{equation}\label{jump-off_cond_dependet}
        \alpha V_{z_1} + \beta W_{z_1} + \gamma=0.
    \end{equation} As expected, $Q_{z_1}=V_{z_{1}}+1$, and $R_{z_1} = W_{z_1}+1$ due to the constraint that $z_1\in W^u(X_-)$.
\end{rem}

We close this section by recapitulating the singular orbit.  $\vp_0$ consists of the five segments below.

\begin{enumerate}[itemsep=7 pt]
    \item [(S1)] $\vp_0$ is the segment of $W^u(X_-)$ between $X_-$ and $z_1$. ($X_-$ is reached as $x\rightarrow-\infty$.) 
    
    \item [(F2)] $\vp_0$ is the heteroclinic orbit for (\ref{fast reduced system}) connecting $(U,P)=(-1,0)$ to $(U,P)=(+1,0)$. The slow coordinates are unchanged.
    
    \item [(S3)] $\vp_0$ is the unique trajectory for (\ref{slow reduced system}) on $M_0^+$ connecting $z_2(=z_1)$ and $z_3(=z_4)$. Note that the $V$ and $W$ coordinates at $z_2$ and $z_3$ are the same, and $(W_{z_3},R_{z_3})=(-W_{z_2},-R_{z_2})$. This ensures both that the jump-off condition (\ref{jump_cond_x}) is met and that $\vp_0$ lands on $W^s(X_-)$ after the fast jump.
    
    \item [(F4)] $\vp_0$ is the heteroclinic orbit for (\ref{fast reduced system}) connecting $(U,P)=(+1,0)$ to $(U,P)=(-1,0)$. The slow coordinates are unchanged.
    
    \item [(S5)] $\vp_0$ is the segment of $W^s(X_-)$ which begins at $q_4$ and approaches $X_-$ as $x\rightarrow\infty$. 
\end{enumerate}
\section{The stability problem and the Maslov index}

Now assume that $(\alpha,\beta,\gamma,D,\eps,x^*)$ satisfy the conditions of Theorem \ref{thm_existence}. Let $\vp_\eps(x)$ be the corresponding standing pulse. We will write $\vp_\eps(x)$ for both the homoclinic orbit of (\ref{standing wave ODE slow}) and the standing pulse solution $(U(x),V(x),W(x))$ of (\ref{PDE model slow}). The goal of this section is to analyze the stability of $\vp_\eps$ by using the Maslov index. We begin by briefly reviewing the stability problem for $\vp_\eps$. Since the following theory is well understood, we omit many of the details. For more background on the stability analysis of solitary waves, we refer the reader to \cite{AGJ,Sandstede02}.

\begin{define}\label{def stability}
    The standing wave $\vp_\eps(x)$ is \textbf{asymptotically stable} relative to (\ref{PDE model slow}) if there is a neighborhood $V\subset BU(\bbR,\bbR^3)$ of $\vp_\eps(x)$ such that if $\psi(x,t)$ solves (\ref{PDE model slow}) with $\psi(x,0)\in V$, then $$||\vp_\eps(x+k)-\psi(x,t)||_\infty\rightarrow 0$$ as $t\rightarrow\infty$ for some $k\in\bbR$.
\end{define}

Although Definition \ref{def stability} refers to nonlinear stability, it is sufficient in this case to prove linear stability \cite{Henry}. Moreover, the essential spectrum is bounded away the imaginary axis in $\bbC=\{z\in\bbC:\mathrm{Re}\,z<0\}$ in our parameter regime \cite[Lemma 3.1]{vHDK_stab}. 

It follows that the stability of $\vp_\eps$ is entirely determined by eigenvalues, which are values $\lambda\in\bbC$ such that 
\begin{equation}\label{eigenvalue_equation_3D}
\begin{aligned}
\lambda u & = \eps^2u_{xx}+u-(3(U(x))^2)u-\eps(\alpha v+\beta w)\\
\lambda \tau v & = v_{xx}+u-v\\
\lambda \theta w & = D^2w_{xx}+u-w
\end{aligned}
\end{equation} has a solution $P_\lambda(x):=(u(x),v(x),w(x))\in BU(\bbR,\bbC^3)$. (Notice that the right hand side of (\ref{eigenvalue_equation_3D}) is the linearization of (\ref{PDE model slow}) about $\vp_\eps(x)$.)  

Due to translation invariance, $\lambda=0$ is necessarily an eigenvalue. However, it is shown in \cite{vHDK_stab} that this eigenvalue is simple. We therefore have the following theorem. \begin{thm}\label{thm stability}
The standing pulse $\vp_\eps(x)$ of (\ref{PDE model slow}) is asymptotically stable if and only if all of the non-zero eigenvalues of (\ref{eigenvalue_equation_3D}) have negative real part.
\end{thm}
To employ the Maslov index, we rewrite (\ref{eigenvalue_equation_3D}) as a first order system on the fast timescale $(\xi)$:
\begin{align}\label{A(lambda,x)}
\dot Y(\xi)=\begin{pmatrix}
u \\ v \\ w \\ p \\ q \\ r
\end{pmatrix}_\xi = 
\begin{pmatrix}
0 & 0 & 0 & 1 & 0 & 0 \\
0 & 0 & 0 & 0 & \epsilon & 0 \\
0 & 0 & 0 & 0 & 0 & \frac{\epsilon}{D} \\
\lambda - 1 + 3U^2 & \alpha \epsilon & \beta \epsilon & 0 & 0 & 0 \\
-\epsilon & \epsilon(\lambda\tau + 1) & 0 & 0 & 0 & 0 \\
-\frac{\epsilon}{D} & 0 & \frac{\epsilon}{D}(\lambda\theta + 1) & 0 & 0 & 0
\end{pmatrix}
\begin{pmatrix}
u \\ v \\ w \\ p \\ q \\ r
\end{pmatrix} =A(\lambda,\xi) Y(\xi).
\end{align}
Notice that we switched the order of the variables in order to make the Hamiltonian structure of the equations clearer.  Conversely, for the nonlinear problem it was more convenient to have the fast and slow variables separated.

We follow the standard approach \cite{AGJ} to analyzing (\ref{A(lambda,x)}), which is to define the stable and unstable bundles of solutions decaying at $-\infty$ and $\infty$ respectively:
\begin{equation}\label{def bundles}
    \begin{aligned}
        E^s(\lambda,\xi) & = \{Y(\xi)\in\bbC^6:Y(\xi) \text{ solves (\ref{A(lambda,x)}) and } Y(\xi)\rightarrow 0 \text{ as } \xi\rightarrow\infty \}\\
        E^u(\lambda,\xi) & = \{Y(\xi)\in\bbC^6:Y(\xi) \text{ solves (\ref{A(lambda,x)}) and } Y(\xi)\rightarrow 0 \text{ as } \xi\rightarrow-\infty \}
    \end{aligned}.
\end{equation} Since the only dependence on $\xi$ of (\ref{A(lambda,x)}) is through $U(\xi)$--which decays exponentially to $-1$ at $\pm\infty$--one would expect the dynamics of (\ref{A(lambda,x)}) to be influenced by those of the constant coefficient system $\dot Y(\xi) = A_\infty(\lambda)Y(\xi)$, where \begin{equation}
    A_\infty(\lambda):=\lim_{\xi\rightarrow\pm\infty}A(\lambda,\xi).
\end{equation} This is indeed the case. Defining $S(\lambda)$ and $U(\lambda)$ to be the stable and unstable subspaces for $A_\infty(\lambda)$ respectively, we have \begin{equation}\label{bundle limits}
    \begin{aligned}
    \lim_{\xi\rightarrow\infty}E^s(\lambda,\xi)&=S(\lambda)\\
    \lim_{\xi\rightarrow-\infty}E^u(\lambda,\xi)&=U(\lambda).
    \end{aligned}
\end{equation}
In particular, the dimensions of $S(\lambda)$ and $U(\lambda)$ are equal to the dimensions of the respective bundles for any $\lambda$ and fixed $\xi$. It is easy to check in this case \cite[\S 3]{vHDK_stab} that $S(\lambda)$ and $U(\lambda)$ are each three-dimensional for all $\lambda$ with non-negative real part. It follows that $E^{u/s}(\lambda,\xi)$ define two-parameter curves in $\Gr_3(\bbC^6)$, or $\Gr_3(\bbR^6)$, if $\lambda\in\bbR$. This fact is the basis for analyzing (\ref{A(lambda,x)}) with the Maslov index in the next section.

\subsection{The Maslov index}
We are now ready to discuss the Maslov index. While we do present all of the necessary definitions and theorems in this paper, the reader may find the more detailed account in \cite[\S 3]{corn} to be helpful. We shall henceforth restrict to $\lambda\in\bbR$. Suppose that $Y_1(\xi)$ and $Y_2(\xi)$ are two solutions of (\ref{A(lambda,x)}) for fixed $\lambda$. Consider the two-form \begin{equation}\label{def omega}
    \omega :=du\wedge dp-\alpha dv\wedge dq -\beta D dw\wedge dr.
\end{equation} It is straightforward (eg., \cite[Theorem 2.1]{corn}) to calculate that \begin{equation}\label{Lambda invariance}
    \frac{d}{d\xi}\omega(Y_1,Y_2)=\omega(Y_1,\dot Y_2)+\omega(\dot Y_1,Y_2) = 0.
\end{equation} In other words, the form $\omega$ is constant in $\xi$ when evaluated on two solutions of (\ref{A(lambda,x)}).\footnote{An ostensibly different symplectic form is used in \cite{corn}. However, this difference is artificial; it is a byproduct of the way that we convert (\ref{eigenvalue_equation_3D}) to a first-order system. (\ref{PDE model slow}) is indeed of the skew-gradient variety studied in \cite{CH,corn}.} Furthermore, if $\alpha$ and $\beta$ are nonzero, $\omega$ is also non-degenerate.  It thus defines a \emph{symplectic form} on $\bbR^6$. Define the symplectic complement of a subspace $V$ as $V^\bot=\{w\in\bbR^6:\omega(v,w)=0,\, \forall v\in V\}$. A (necessarily three-dimensional) subspace $V\subset\bbR^6$ such that $V=V^\bot$ is called \emph{Lagrangian}. We can therefore rephrase (\ref{Lambda invariance}) as follows: The set of Lagrangian planes, denoted $\Lambda(3)$, is an invariant set of the dynamical system induced on $\Gr_3(\bbR^6)$ by (\ref{A(lambda,x)}).

The set $\Lambda(3)$ happens to be a submanifold of $\Gr_3(\bbR^6)$ of dimension 6. The interesting fact for our purposes is that $\pi_1(\Lambda(3))=\bbZ$, which means that an integer winding number can be associated with curves in this space. (By contrast, $\pi_1(\Gr_3(\bbR^6))=\bbZ/2\bbZ$, so no such integer index exists in the full Grassmannian.) This winding number is the Maslov index.

Since we will be dealing with non-closed curves in general, we will define the Maslov index as an intersection number instead of a winding number (\`{a} la Poincar\'{e} duality \cite{hatcher}). More precisely, we will fix a Lagrangian plane $V$--called the reference plane--and count the number of times that a curve of Lagrangian planes intersects $V$. Arnol'd established this definition of the Maslov index in the case of one-dimensional intersections \cite{Arnold67}, and then Robbin and Salamon generalized it for intersections of any dimension \cite{RS93}.

To handle the accounting for multidimensional intersections, Robbin and Salamon developed the ``crossing form." This is a quadratic form whose signature determines the contribution to the Maslov index at each point. We remind the reader that the signature of a quadratic form $Q$ is the difference of the positive and negative indices of inertia: \begin{equation}
    \mathrm{sign}(Q)=n_+(Q)-n_-(Q).
\end{equation}
For a reference plane $V\in\Lambda(3)$ and a curve of Lagrangian planes $\gamma(t)\in \Lambda(3)$, a \emph{conjugate point} is a time $t_*$ such that $\gamma(t_*)\cap V\neq \{0\}$. The crossing form is defined on the intersection, so the contribution to the Maslov index at a $k-$dimensional crossing can be anywhere between $-k$ and $k$. A crossing is called \emph{regular} if $Q$ is non-degenerate.

Below we will formally define the Maslov index of the solitary wave $\vp_\eps$. The curve of interest is the unstable bundle $E^u(0,\lambda)$, which we can think of as containing all potential 0-eigenfunctions (i.e., those which satisfy the left boundary conditions). In the spirit of a shooting argument, the reference plane should be the right boundary data for a potential 0-eigenvector, namely $S(0)$; see (\ref{bundle limits}). Actually, for technical reasons, it is untenable to consider $E^u(0,\xi)$ on all of $\bbR$ due to the presence of a conjugate point at $+\infty$. (See Remark \ref{rem endpoints}.) The rigorous definition of the Maslov index for a standing wave--due to Chen and Hu \cite{CH}--fixes this issue by truncating the curve and using the stable bundle at the new endpoint as the reference plane. 
\begin{define}\label{def Maslov of wave}
    Let $\xi_\infty$ be large enough so that \begin{equation}\label{cutoff} U(0)\cap E^s(0,\xi)=\{0\}\text{ for all } \xi\geq \xi_\infty.\end{equation} We define the \textbf{Maslov index} of $\vp_\eps$ to be \begin{equation}
        \Maslov:=\sum_{\xi_*\in(-\infty,\xi_\infty)}\mathrm{sign}\,\Gamma(E^u,E^s(0,\xi_\infty),\xi_*)+n_+(\Gamma(E^u,E^s(0,\xi_\infty),\xi_\infty)),
    \end{equation} where the sum is taken over all interior crossings of $\xi\mapsto E^u(0,\xi)$ with $\Sigma$, the train of $E^s(0,\xi_\infty)$. If $\xi_*$ is a conjugate point, and $\psi\in E^u(0,\xi_*)\cap E^s(0,\xi_\infty)$, then the crossing form is given by  \begin{equation}\label{def xing form}
        \Gamma(E^u,E^s(0,\xi_\infty),\xi_*)=\omega(\psi,A(0,\xi_*)\psi).
    \end{equation}
\end{define}

\begin{rem}\label{rem endpoints}
    We only count negative crossings at $-\infty$ (of which there are none) and positive crossings at $\xi_\infty$ by convention. Notice that there is a guaranteed to be a crossing at $\xi=\xi_\infty$, since $\vp'_\eps(\xi_\infty)\in E^s(0,\xi_\infty)\cap E^u(0,\xi_\infty)$. This crossing is one-dimensional--equivalently, the translation eigenvalue is simple--by virtue of the transverse construction of the wave.
\end{rem}

A key feature of Definition \ref{def Maslov of wave} is that $\lambda=0$ is fixed. One can typically show that $\Maslov$ either counts real, unstable eigenvalues or gives a lower bound on the count. Moreover, (\ref{A(lambda,x)}) is the equation of variations for (\ref{standing wave ODE fast}) when $\lambda=0$, which means that one can glean spectral information from how the wave itself is constructed. Indeed, this is the premise for the calculation in the next section.

The next step would be to prove that the Maslov index actually counts all unstable eigenvalues. Since our focus here is on the calculation of the index more than stability per se, we will not go down that path. However, the authors in \cite{vHCNT} verified that the Maslov index does indeed give the desired count. For a blueprint on how to prove this equality in general, we refer the reader to \cite{CH14,corn}.
\section{Calculation of the Maslov index}

We are now prepared to calculate the Maslov index using the framework developed in \cite{CJ18}. In order to motivate the calculation, we first state the known stability result for $\vp_\eps$.
\begin{thm}[Theorem 4.1 of \cite{vHDK_stab}, Theorem 1 of \cite{vHCNT}]\label{thm stability result}
    Let $\vp_\eps(\xi)$ be a standing single pulse solution of (\ref{PDE model slow}), as described in Theorem \ref{thm_existence}.  Then $\vp_\eps$ is stable in the sense of Definition \ref{def stability} if and only if \begin{equation}\label{stab crit}
    \alpha V_0 + \frac{\beta}{D} W_0<0,
    \end{equation} where $V_0$ and $W_0$ are the values (to leading order in $\eps$) of $V,W$ at the jump-off points $z_i$.
\end{thm}
\begin{rem}\label{rem stab1}
    Recall that $V_0,W_0\in(-1,0)$, so at least one of $\alpha,\beta$ must be negative in order for $\vp_\eps$ to be unstable.
\end{rem}
\begin{rem}
    Values of $(\alpha,\beta,D)$ yielding equality in (\ref{stab crit}) correspond to a saddle-node bifurcation of homoclinic orbits \cite[Theorem 2.1]{vHDK_stab}.
\end{rem}

Our perspective in this section will be to consider $\alpha$ and $\beta$ as parameters. In light of Theorem \ref{thm stability result}, we expect to see a conjugate point appear or disappear as $\alpha$ and $\beta$ are tuned to cause a sign change in (\ref{stab crit}).

The ensuing calculation rests on two ideas.  The first is that the stable and unstable bundles (\ref{def bundles}) are everywhere tangent to the stable and unstable manifolds for $X_\eps^-$ as a fixed point of (\ref{standing wave ODE fast}). The second is that the nonlinear objects $W^{u/s}(X_\eps^-)$ can be tracked using Fenichel theory \cite{Fen79,Jones_GSP}. We know that $E^u(0,\xi)$ is three-dimensional, with one direction given by the wave velocity $\vp'_\eps$. Using GSPT, we are able to discern the other two dimensions using the geometry of phase space as opposed to having to solve the non-autonomous linear system (\ref{A(lambda,x)}).

The crux of GSPT is that the critical manifolds $M_0^{\pm}$ perturb to locally invariant manifolds $M_\eps^{\pm}$ when $0<\eps\ll1$. Moreover, the flow on $M_\eps^{\pm}$, to leading order, is is given by (\ref{slow reduced system}). Actually, something even stronger is true.  Recall that $M_0^{\pm}$ is the union of fixed points of the fast system (\ref{fast reduced system}). Each of these fixed points has one-dimensional stable and unstable manifolds, obtained by linearizing (\ref{fast reduced system}) at $U=\pm1$. It is therefore possible to define $W^{u/s}(M_0^{\pm})$ as the union of the corresponding invariant manifolds for the points comprising $M_0^{\pm}$. These, too, perturb to locally invariant manifolds $W^{u/s}(M_\eps^{\pm})$. The fact that (\ref{fast reduced system}) is 2D Hamiltonian and (\ref{slow reduced system}) is essentially linear will make it quite easy to describe these objects.

The smooth convergence of the slow manifolds and their attendant invariant manifolds allows us to work primarily with $\eps=0$, provided that all conjugate points are regular. However, the transitions from fast-to-slow dynamics and vice-versa require some care. Indeed, the Maslov index calculation is properly happening in the tangent bundle along the wave. Although the wave itself is continuous in the limit, its derivative (and more generally the tangent space to $W^u(X_0^{-})$) has a jump discontinuity at these transitions. To figure out the $\eps>0$ links in $\Lambda(3)$ between the fast and slow pieces at $z_i$, it is shown in \cite{CJ18} that one must treat (\ref{standing wave ODE fast}) as a constant coefficient system at $z_i$. As we shall see, these `corners' are the most interesting part of the calculation since the bifurcation manifests itself there.

In light of the preceding paragraph, there should be nine segments (= five components of the wave and four corners) that we must analyze to calculate the index. However, our flexibility in choosing the right endpoint of the unstable bundle will shorten the calculation considerably. Since a similar calculation is done in full detail in \cite{CJ18}, we will omit many of the details and focus on the pieces that influence stability. In referring to the corners, we will use the notation ``$i.5$" for the transition from segment $i$ to $i+1$.

\subsection{The reference plane}\label{sec ref plane}
To utilize Definition \ref{def Maslov of wave}, we need to select a reference plane. Selecting a reference plane is tantamount to choosing a value $\xi_\infty$ so that $U(0)\cap E^s(0,\xi)=\{0\}$ for all $\xi\geq \xi_\infty$. To that end, we begin by linearizing (\ref{standing wave ODE fast}) about $X_-$. We can obtain the leading-order terms by setting $\eps=0$ separately in (\ref{fast reduced system}) and (\ref{slow reduced system}). This way the $(U,P),(V,Q)$, and $(W,R)$ decouple nicely, and we we compute the following eigenvectors in $(u,v,w,p,q,r)$ coordinates: \begin{equation}\label{evecs_minus1}
    S(0) = \mathrm{span}\left\{\begin{pmatrix}
    1\\0\\0\\-\sqrt{2}\\0\\0
    \end{pmatrix}, \begin{pmatrix}
    0\\1\\0\\0\\-1\\0
    \end{pmatrix},\begin{pmatrix}
    0\\0\\1\\0\\0\\-1
    \end{pmatrix}\right\},\hspace{5pt} U(0) = \mathrm{span}\left\{\begin{pmatrix}
    0\\0\\1\\0\\0\\1
    \end{pmatrix},\begin{pmatrix}
    0\\1\\0\\0\\1\\0
    \end{pmatrix}, \begin{pmatrix}
    1\\0\\0\\\sqrt{2}\\0\\0
    \end{pmatrix}\right\}.
\end{equation}

Without loss of generality, we assume that $D>1$. With this assumption, the eigenvectors (\ref{evecs_minus1}) appear from left to right in order of increasing eigenvalue. We label these $\eta_1, \dots, \eta_6$, with corresponding eigenvalues $\mu_i$. To leading order in $\eps$, the eigenvalues are \begin{equation}
    \begin{aligned}
    \mu_1 &= -\sqrt{2}\\
    \mu_2 &= -\eps\\
    \mu_3 &= -\eps/D\\
    \mu_4 &= \eps/D\\
    \mu_5 &= \eps\\
    \mu_6 &= \sqrt{2}\\
    \end{aligned}
\end{equation}

To choose the reference plane, we flow $E^s(0,\xi)$ backwards in time from $+\infty$ until we reach a convenient point. Because (\ref{slow reduced system}) itself is affine linear, there is no change to the slow directions $\eta_2$ and $\eta_3$ along segment $(S5)$. Neither are there any changes along the fast flow, so we know that $\eta_2,\eta_3\in E^s(0,\xi)$ for any point along the back. The fast stable direction also does not change along $(S5)$; it is given by $\eta_1$, which is the direction in which $\vp_0$ approaches $z_4$ along the back. This proves that $E^s(0,\xi)$ is tangent to $W^s(M_0^-)$ at the landing point $z_4$. Along the back itself, the stable direction is tangent to the heteroclinic orbit. This can be computed from (\ref{up_hetero}) as $(1, \sqrt{2}U)$ in $(u,p)$ space.

Let $\xi$ be some time for which $\vp_\eps(\xi)$ is on the back. To verify condition (\ref{cutoff}), we compute the determinant
\begin{equation}
    \det[U(0),E^s(0,\xi)]=\det\begin{pmatrix}
    0 & 0 & 1 & 1 & 0 & 0\\
    0 & 1 & 0 & 0 & 1 & 0\\
    1 & 0 & 0 & 0 & 0 & 1\\
    0 & 0 & \sqrt{2} & \sqrt{2}U & 0 & 0\\
    0 & 1 & 0 & 0 & -1 & 0\\
    1 & 0 & 0 & 0 & 0 & -1\\
    \end{pmatrix}=4\sqrt{2}(1-U).
\end{equation}
Evidently we are free to select any point along the back besides $\pm1$, so we will pick the jump midpoint where $U=0$. We will therefore count intersections of $E^u(0,\xi)$ with the reference plane \begin{equation}\label{def ref plane}
    V=\mathrm{span}\left\{\begin{pmatrix}
    1\\0\\0\\0\\0\\0
    \end{pmatrix}, \begin{pmatrix}
    0\\1\\0\\0\\-1\\0
    \end{pmatrix},\begin{pmatrix}
    0\\0\\1\\0\\0\\-1
    \end{pmatrix}\right\}=\mathrm{span}\{\eta_1+\eta_6,\eta_2,\eta_3\}.
\end{equation}

\subsection{The fast front}\label{calc fast front}
This sections begins the calculation of the Maslov index. As mentioned earlier, $E^u(0,\xi)$ is everywhere tangent to $W^u(X_\eps^-)$ along $\vp_\eps(\xi)$. By the same reasoning as \S \ref{sec ref plane}, the tangent space to $W^u(X_\eps^-)$ is spanned by $\{\eta_4,\eta_5,\eta_6\}$ throughout segment $(S1)$ and corner $1.5$. It is clear from (\ref{def ref plane}) that there are no conjugate points for these pieces.

Along the fast front, the slow directions remain the same, and the fast unstable direction is given by the tangent vector $\vp_\eps$. Again using (\ref{up_hetero}), we find the $(u,p)$ components to be $(1,-\sqrt{2}U)$. To detect conjugate points (i.e., intersections with $V$), we compute \begin{equation}
    \det[E^u(0,\xi),V]=-4\sqrt{2}U,
\end{equation} which is zero if and only if $U=0$. This is the midpoint of the jump, just like we chose for the reference plane. If we selected a different $\xi_\infty$ for Definition \ref{def Maslov of wave}, then this conjugate point would have moved accordingly. To compute the dimension of this crossing, we first verify that the intersection is one-dimensional, spanned by \begin{equation}
    \psi=\begin{pmatrix}
    1\\0\\0\\0\\0\\0
    \end{pmatrix}.
\end{equation} We next apply (\ref{def xing form}) with $\eps=0$ and $\xi$ such that $U=0$ to see that \begin{equation}
    \omega(\psi,A(0,\xi)\psi)=du\wedge dp(\psi,A(0,\xi)\psi)=-1.
\end{equation} Thus the crossing form is one-dimensional and negative definite, so the contribution to the Maslov index along the fast front $(F2)$ is $-1.$

\subsection{Passage near the right slow manifold}\label{calc slow mfd}

It turns out that the stability result will hinge on corner $2.5$, so we save that section until the end. We now consider the passage of $E^u(0,\xi)$ near $M_\eps^+$. The slow flow takes over for this section, so the fast direction will be constant to leading order. As one would expect, it is the unstable fast direction $\eta_6$ that persists for this segment; see \cite[\S 4.2-4.3]{CJ18}.

The slow directions are slightly more complicated. Recall the jump-off condition (\ref{jump_cond_x}) for the fast jump. This condition describes a one-dimensional set of valid initial conditions on $M_\eps^+$. By flowing this set forward in time, we obtain a two-dimensional manifold with boundary $S$ foliated by the slow trajectories. The tangent space to this manifold along $(S3)$ gives the slow directions of $T_{\vp_\eps(\xi)}W^u(X_\eps^-)$ near this segment.

To get an explicit expression, we first solve (\ref{slow reduced system}) generically on $M_0^+$ (i.e., with $U\equiv1$):
\begin{equation}\label{slow flow solution}
    \begin{aligned}
    V(x) & = c_1 e^x + c_2e^{-x} + 1 \\
Q(x) & = c_1 e^x - c_2e^{-x} \\
W(x) & = c_3 e^{x/D} + c_4 e^{-x/D} + 1 \\
R(x) & = c_3 e^{x/D} - c_4 e^{-x/D}
    \end{aligned}
\end{equation} Because the steady state equation for (\ref{PDE model slow}) has reversibility symmetry, we know that $V$ and $W$ will achieve their maxima at the same point. We may assume that this maximum is $x=0$ for (\ref{slow flow solution}), from which we conclude that $c_1=c_2$ and $c_3=c_4$. The solutions of interest for (\ref{slow flow solution}) satisfy the equations \begin{equation}
    \begin{aligned}
    (V-1)^2-Q^2&=4c_1^2\\
    (W-1)^2-R^2&=4c_3^2,
    \end{aligned}
\end{equation} which intersect the lines $Q=V+1$ and $R=W+1$ from (\ref{jump_off_coords}) when 
\begin{equation}\label{VW_connection_slow_mfd}
    \begin{aligned}
    V&=-c_1^2\\
    W&=-c_3^2.
    \end{aligned}
\end{equation}
We combine (\ref{VW_connection_slow_mfd}) with (\ref{jump-off_cond_dependet}) to see that \begin{equation}\label{c_i jumpoff cond}
    \alpha c_1^2+\beta c_3^2=\gamma,
\end{equation} which we can solve for $c_1$, knowing that both $c_{1,3}<0$.
We can therefore describe $S$ as the graph of a function of $x$ and $c_1$: \begin{equation}
    (V,Q,W,R)=h(x,c_1)=\begin{pmatrix}
c_1\cosh(x) + 1 \\
c_1\sinh(x) \\
c_3(c_1)\cosh\left(\frac{x}{D}\right) \\
c_3(c_1)\sinh\left(\frac{x}{D}\right)
\end{pmatrix}
\end{equation} We differentiate to determine a basis for the tangent space to $S$: \begin{equation}\label{S_tangent_space}
T_{\vp_0(x)} S= \mathrm{span}\left\{\frac{\p h}{\p x},\frac{\p h}{\p c_1}\right\}=\mathrm{span}\begin{pmatrix}
c_1\sinh(x) & \cosh(x)\\
c_1\cosh(x) & \sinh(x)\\
\frac{c_3}{D}\sinh\left(\frac{x}{D}\right) & c_3'\cosh\left(\frac{x}{D}\right)\\
\frac{c_3}{D}\cosh\left(\frac{x}{D}\right) & c_3'\sinh\left(\frac{x}{D}\right)
\end{pmatrix}
\end{equation} where $'=\frac{d}{dc_1}$. Since the fast direction is unchanged, it is clear that we can check for conjugate points by seeing if (\ref{S_tangent_space}) intersects the plane spanned by $\eta_2$ and $\eta_3$. To that end, a somewhat messy computation shows that  \begin{equation}\label{slow_mfd_det}
    \det\left[\eta_2,\eta_3,\frac{\p h}{\p x},\frac{\p h}{\p c_1}\right]=(c_3/D-c_1c_3')\,e^{x(1+1/D)}.
\end{equation} Clearly, the sign of (\ref{slow_mfd_det}) is independent of $x$. Differentiating (\ref{c_i jumpoff cond}) implicitly, we find that \begin{equation}\label{slow_mfd_det_2}
    \frac{c_3}{D}-c_1c_3'=\frac{c_3}{D}+\frac{\alpha c_1^2}{\beta c_3}=\frac{(\beta/D)c_3^2+\alpha c_1^2}{\beta c_3}.
\end{equation} We must now evaluate this expression along $\vp_0(x)$. Using (\ref{jump-off_cond_dependet}) and (\ref{VW_connection_slow_mfd}), the numerator of (\ref{slow_mfd_det_2}) simplifies to \begin{equation}
    -\left(\alpha V_{z_1}+\frac{\beta}{D}W_{z_1}\right).
\end{equation} Notice that this is exactly the expression appearing in the stability condition (\ref{stab crit}). We are assuming that we are not in the regime where the bifurcation occurs, so (\ref{slow_mfd_det}) does not vanish anywhere along $(S3)$. Hence there are no conjugate points, and the Maslov index is unchanged along this segment. Although this is true regardless of the sign, one should be suspicious of the stability criterion appearing in this way.

\subsection{The fast back}

Corner 3.5 turns out to be trivial after the pain of the previous section.  The fast direction is already in the correct position for the fast back--tangent to $W^u(X_+)$--from the passage near $M_\eps^+$. An application of the $(k+\sigma)$ Exchange Lemma \cite[\S 6.1]{Jones_GSP} with $\sigma=2$ shows that the slow directions also remain unchanged. We therefore have no contribution to the Maslov index.

For the back, we can focus entirely on the two-dimensional fast system.  Indeed, the calculations (\ref{slow_mfd_det}) and (\ref{slow_mfd_det_2}) show that the slow directions at launch (which don't change over the back) are transverse to the slow directions in $V$.

We again use (\ref{up_hetero}) to get a $U$-dependent expression for the fast direction of $W^u(X_+)$. This time, we use the branch of (\ref{up_hetero}) with negative $P$, so the tangent vector is given by $(u,p)=(1,\sqrt{2}P)$. As in \S \ref{calc fast front}, there is exactly one conjugate point at $P=0$. We could easily compute the crossing form again, but we'll instead argue geometrically that the crossing is negative. The two fast jumps form a separatrix for (\ref{fast reduced system}). For a planar system, the Maslov index measures the winding of the angle that the tangent vector to the wave makes. It is clear that the tangent vector is rotating clockwise at both of $(0,1/\pm\sqrt{2})$, so the crossing must have the same direction for both conjugate points.

Although the crossing is regular and negative, this conjugate point does not contribute to the Maslov index because we selected $\xi_\infty$ so that this is an endpoint crossing. By Definition \ref{def Maslov of wave}, only positive crossings are counted at the right endpoint. 

Excluding corner 2.5, it follows that $\Maslov=-1$, with the only contributing conjugate point being the midpoint of the fast front. In the final section, we analyze corner 2.5.  We now know that in order for the wave to be stable, there must be a conjugate point in the positive direction to offset the one from the front.

\subsection{Corner 2.5: Arrival at the right slow manifold}\label{sec corner}

As the wave $\vp_0(x)$ approaches $z_2$, we know from \S \ref{calc fast front} that $W^u(X_-)$ will be tangent to \begin{equation}
    Y=\mathrm{span}\{\eta_1,\eta_4,\eta_5\}.
\end{equation} The stable direction $\eta_1$ is present because it is the limit as $U\rightarrow 1$ of $(1,-\sqrt{2}U)$. In other words, the fast front much approach the fixed point $(1,0)$ tangent to the stable manifold of the fixed point.

On the other hand, the jump-off condition (\ref{jump-off_cond_dependet}) must be satisfied by incoming points if they are to spend a long time near $M_\eps^+$. This forces an abrupt reorientation of $W^u(X_-)$ right at $z_2$. The fast direction will be unstable, and the slow directions are given by (\ref{S_tangent_space}) evaluated at the landing point. It is not difficult to see that we can write these in terms of the $\eta_i$ as \begin{equation}\label{def Z}
    Z = \mathrm{span}\left\{ \eta_6, \eta_2 + \frac{1}{D} \eta_3 + \frac{W_{z_1}}{D}  \eta_4 + V_{z_1} \eta_5, -\frac{\alpha}{\beta} \eta_4 + \eta_5 \right\}.
\end{equation}

Our task is to figure out the $\eps>0$ smooth connection between $Y$ and $Z$ in $\Lambda(3)$. To do this, we first observe that the tangent space to $W^u(X_-^\eps)$ evolves on the fast timescale because it always contains at least one fast direction. It follows that this reorientation must happen arbitrarily close to $z_1$ as $\eps\rightarrow0$. The relevant dynamical system to study is therefore \begin{equation}\label{linearization_corner}
    Y'(x)=BY(x),
\end{equation} where $B=A(0,\xi_*)$, and $U(\xi_*)=1$. System (\ref{linearization_corner}) induces an equation on $\Gr_3(\bbR^6)$, which we know from (\ref{Lambda invariance}) has $\Lambda(3)$ as an invariant submanifold. The phase portrait of (\ref{linearization_corner}) is fairly simple and described in detail in \cite[Appendix B]{CJ18}. The highlights are that the fixed points are direct sums of eigenspaces, each of which is hyperbolic. This includes $Y$. We will find the orbit connecting $Y$ and $Z$ by treating (\ref{linearization_corner}) as a boundary value problem. In this sense, our analysis is similar in spirit to the Brunovsky approach to inclination lemmas \cite{Bru96,Bru99}.

Since $Y$ is a fixed point, any candidate link between $Y$ and $Z$ must lie in $W^u(Y)\subset\Lambda(3)$. Using row reduction, we can compute that all such planes must be expressible in the form $\mathrm{span}\{v_1,v_2,v_3\}$, with \begin{equation}\label{unstable mfd of Y}
\begin{aligned}
v_1 & = \eta_1+a_{12}\eta_2+a_{13}\eta_3+a_{16}\eta_6\\
v_2 & = \eta_4+a_{46}\eta_6\\
v_3 & = \eta_5+a_{56}\eta_6.
\end{aligned}
\end{equation} Actually, (\ref{unstable mfd of Y}) gives $W^u(Y)$ in the full Grassmannian.  We can reduce this set to three dimensions by applying the form $\omega$ pairwise on the $v_i$ and imposing the condition that it vanishes. As a result, we can express $W^u(Y)\subset\Lambda(3)$ as \begin{equation}\label{Y unstable mfd Lag}
W^u(Y)=\mathrm{sp}\{\eta_1+\delta_2\eta_2+\delta_3\eta_3+\delta_6\eta_6,\eta_4+\frac{\beta D}{\sqrt{2}}\delta_3\eta_6,\eta_5+\frac{\alpha}{\sqrt{2}}\delta_2\eta_6\},
\end{equation} for $\delta_{2,3,6}\in\mathbb{R}$.

Comparing (\ref{Y unstable mfd Lag}) and (\ref{def Z}), it appears that we are in trouble since $Z\notin W^u(Y)$. However, these expressions all hold to leading order only, and we shall see that there are 3-planes arbitrarily close to $Z$ which do lie in $W^u(Y)$. To prove this, we express each of $Z$ and a generic plane in $W^u(Y)$ in Pl\"{u}cker coordinates \cite[Appendix A]{CJ18} using the basis $\{\eta_i\}$:

\begin{center}
\renewcommand*{\arraystretch}{2}
\begin{longtable}{|c|c|c|}
\hline
Basis Vectors & $W^u(Y)$ & $Z$ \\
\endfirsthead
\hline
Basis Vectors & $W^u(Y)$ & $Z$ \\
\hline
\endhead
\hline
$\{1,2,3\}$ & 0 & 0 \\
\hline
$\{1,2,4\}$ & 0 & 0 \\
\hline
$\{1,2,5\}$ & 0 & 0 \\
\hline
$\{1,2,6\}$ & 0 & 0 \\
\hline
$\{1,3,4\}$ & 0 & 0 \\
\hline
$\{1,3,5\}$ & 0 & 0 \\
\hline
$\{1,3,6\}$ & 0 & 0 \\
\hline
$\{1,4,5\}$ & 1 & 0 \\
\hline
$\{1,4,6\}$ & $\displaystyle\frac{\alpha}{\sqrt{2}} \delta_2$ & 0 \\
\hline
$\{1,5,6\}$ & $-\displaystyle\frac{\beta D}{\sqrt{2}} \delta_3$ & 0 \\
\hline
$\{2,3,4\}$ & 0 & 0 \\
\hline
$\{2,3,5\}$ & 0 & 0 \\
\hline
$\{2,3,6\}$ & 0 & 0 \\
\hline
$\{2,4,5\}$ & $\delta_2$ & 0 \\
\hline
$\{2,4,6\}$ & $\displaystyle\frac{\alpha}{\sqrt{2}} \delta_2^2$ & $-\displaystyle\frac{\alpha}{\beta}$ \\
\hline
$\{2,5,6\}$ & $-\displaystyle\frac{\beta D}{\sqrt{2}} \delta_2 \delta_3$ & 1 \\
\hline
$\{3,4,5\}$ & $\delta_3$ & 0 \\
\hline
$\{3,4,6\}$ & $\displaystyle\frac{\alpha}{\sqrt{2}} \delta_2 \delta_3$ &$-\displaystyle\frac{\alpha}{\beta D}$ \\
\hline
$\{3,5,6\}$ & $-\displaystyle\frac{\beta D}{\sqrt{2}} \displaystyle\delta_3^2$ & $\displaystyle\frac{1}{D}$ \\
\hline
$\{4,5,6\}$ & $\delta_6$ & $\displaystyle\frac{\beta(W_{z_1}/D)+\alpha V_{z_1}}{\beta}$ 
\\ \hline
\end{longtable}
\end{center}
Pl\"{u}cker coordinates are projective, so only the ratio of terms has meaning. We use the notation $p_{ijk}$ for the Pl\"{u}cker coordinates of the plane spanned by $\eta_i,\eta_j,$ and $\eta_k$. The goal is to choose $\delta_{2,3,6}$ so that the resulting Lagrangian 3-plane is arbitrarily close to $Z$. We first compute the ratio \begin{equation}\label{delta3}\frac{p_{246}}{p_{346}}=\frac{\delta_2}{\delta_3}=\frac{-\alpha/\beta}{-\alpha/(\beta D)} =D.\end{equation} One can verify that this choice is consistent with any other $p_{ijk}$ that are degree 2 in $\delta_2$ and $\delta_3$. To fix $\delta_6$, we compute the ratio \begin{equation}\frac{p_{246}}{p_{456}}=\frac{\alpha\delta_2^2}{\sqrt{2}\delta_6}=\frac{-\alpha}{\beta(W_{z_1}/D)+\alpha V_{z_1}},\end{equation} from which we conclude that \begin{equation}\label{delta6}
    \delta_6 = -\frac{\beta(W_{z_1}/D)+\alpha V_{z_1}}{\sqrt{2}}\delta_2^2.
\end{equation}
Choosing $\delta_3$ and $\delta_6$ to satisfy (\ref{delta3}) and (\ref{delta6}) ensures that the Pl\"{u}cker coordinates have the proper ratios whenever the corresponding coordinates for $Z$ are nonzero. Notice that the remaining non-zero coordinates for the generic 3-plane in $W^u(Y)$ are all first order in $\delta_2$ or $\delta_3$. It follows that we can make the \emph{projective} coordinates as close to 0 as desired by sending $\delta_2$ towards infinity. This completes the proof.

Now that we have a plane $X\in W^u(Y)$ that is arbitrarily close to $Z$, all that remains is to compute the trajectory (in backwards time) through $X$. This is easily done since the $\eta_i$ are eigenvectors of $B$. Denoting the trajectory through (\ref{Y unstable mfd Lag}) as $\Phi(x)$, we see that 

\begin{equation}\label{trajectory}
\Phi(x)=\mathrm{span}\left\{
\begin{array}{c}
e^{\mu_1 x} \eta_1 + \delta_2 e^{\mu_2 x} \eta_2 + \delta_3 e^{\mu_3 x} \eta_3 + \delta_6 e^{\mu_6 x} \eta_6, \\
e^{\mu_4 x} \eta_4 + \frac{\beta D}{\sqrt{2}} \delta_3 e^{\mu_6 x} \eta_6, \\
e^{\mu_5 x} \eta_5 + \frac{\alpha}{\sqrt{2}} \delta_2 e^{\mu_6 x} \eta_6
\end{array}
\right\},
\end{equation}
for $x\in(-\infty,0]$. To check for conjugate points, we evaluate the determinant of $[V, \Phi(x)],$ taking care to write $V$ in the basis $\{\eta_i\}$: \begin{equation}\label{corner_conj}
    \det\begin{pmatrix}
1 & 0 & 0 & 0 & 0 & e^{\mu_1 x} \\
0 & 1 & 0 & 0 & 0 & \delta_2 e^{\mu_2 x} \\
0 & 0 & 1 & 0 & 0 & \delta_3 e^{\mu_3 x} \\
0 & 0 & 0 & e^{\mu_4 x} & 0 & 0 \\
0 & 0 & 0 & 0 & e^{\mu_5 x} & 0 \\
1 & 0 & 0 & \frac{\beta D}{\sqrt{2}} \delta_3 e^{\mu_6 x} & \frac{\alpha}{\sqrt{2}} \delta_2 e^{\mu_6 x} & \delta_6 e^{\mu_6 x}
\end{pmatrix} = e^{(\mu_4+\mu_5)x}(\delta_6 e^{\mu_6x}-e^{\mu_1x}).
\end{equation}

Using the fact that $\mu_1=-\mu_6$ (for any $\eps$), we see that (\ref{corner_conj}) vanishes if \begin{equation}\label{conj condition}
\delta_6=e^{2\mu_1x}.
\end{equation} The right hand side of (\ref{conj condition}) takes all values in $[1,\infty)$. Since $|\delta_6|\gg 1$ to be close to $Z$, it follows that there is a conjugate point if and only if $\delta_6>0$, or equivalently, \begin{equation}\label{conj_pt_crit}
    \beta(W_{z_1}/D)+\alpha V_{z_1}<0.
\end{equation}

Comparing Theorem \ref{thm stability result} and (\ref{conj_pt_crit}), the conclusion is that there is a conjugate point at this corner if and only if the wave is stable. We know from the previous sections that $\Maslov=-1$ excluding this corner, so it must be the case that this conjugate point crosses in the positive direction in order to get $\Maslov=0$ for the whole wave.

To verify this, we observe that the intersection of $V$ and $\Phi(x)$ is spanned by \begin{equation}
    \psi = e^{\mu_1x}\eta_1+\delta_2e^{\mu_2x}\eta_2+\delta_3e^{\mu_3x}\eta_3+e^{\mu_1x}\eta_6.
\end{equation} At this corner, $A(0,\xi)=B$, and $B\eta_i=\mu_i\eta_i$. It is straightforward then to calculate that \begin{equation}\label{corner xing form}
    \begin{aligned}
    \omega(\psi,A(0,\xi)\psi) & = \omega(e^{\mu_1x}\eta_1+e^{\mu_1x}\eta_6, \mu_1e^{\mu_1x}\eta_1+\mu_6e^{\mu_1x}\eta_6)\\
    &= e^{2\mu_1x}(\mu_6-\mu_1)\omega(\eta_1,\eta_6)\\
    & = 2\sqrt{2}e^{2\mu_1x}(\mu_6-\mu_1)>0.
    \end{aligned}
\end{equation}
\section{Conclusion}

The calculation of the previous section sheds new light on how the stability criterion (\ref{stab crit}) manifests itself in the construction of the wave. Considering (\ref{slow_mfd_det_2}) and (\ref{Y unstable mfd Lag}), the sign of $\delta_6$ (i.e., the stability criterion) is telling us about the exchange of the fast directions near $M_\eps^+$. To better understand this, we consider the projection of $E^u(0,\xi)$ onto the fast subsystem (\ref{fast reduced system}).

\begin{figure}
\centering
\includegraphics[scale = 0.6]{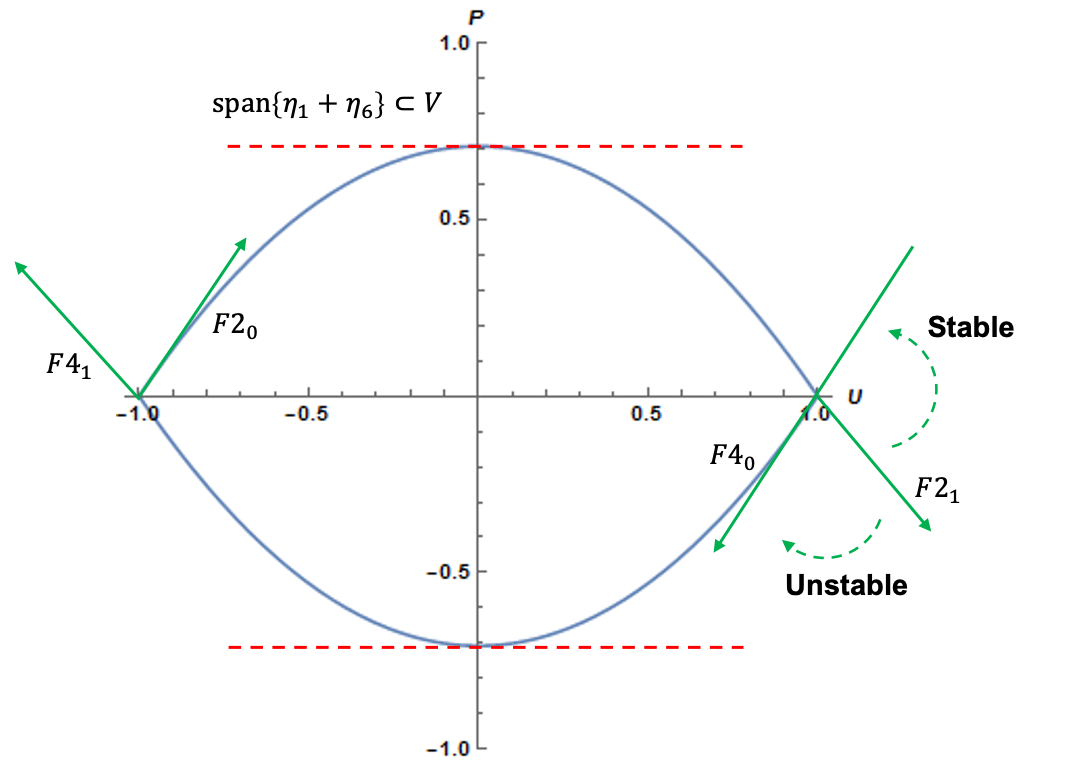}
\caption{The Maslov index is a winding number for $\vp'_0(\xi)$ in the $UP-$plane. The inequality (\ref{stab crit}) reflects the two ways in which $\vp'_0$ can reorient itself for the back at $(U,P)=(1,0)$.}
\label{fig_maslov_planar}
\end{figure}

The fast component of $E^u(0,\xi)$ along the jumps is tangent to $\vp_0$. Along both $F2$ and $F4$ the angle of the tangent line rotates clockwise as $\xi$ increases. As a result, the crossing form is negative definite at all crossings. The challenge is what happens at $(U,P)=(1,0)$, where the tangent vector at the end of $F2$ and the beginning of $F4$ are different. See Figure \ref{fig_maslov_planar}.

At the endpoint of $F2$ ($F2_1$ in the figure), $\vp_0'(\xi)$ is tangent to $\eta_1$. In order to get back to $\eta_6$ for segment $F4$, it can either continue rotating clockwise or reverse course and go counter-clockwise. In the latter case, the net angular rotation along the front ends up being 0, since the fixed points $U=\pm1$ have the same eigenvectors. Moreover, the tangent vector must necessarily be horizontal again (i.e., parallel to $\eta_1+\eta_6$), and it must cross the horizontal subspace in the positive direction. This is exactly what we observed in \S \ref{sec corner}; although the intersection of $E^u(0,\xi)$ with $V$ had slow components, the crossing form calculation (\ref{corner xing form}) reduced to a calculation in the fast components. As expected, the crossing form was positive definite, so the winding along the front was ``undone" at the corner. Conversely, in the unstable case we observed no conjugate points at corner 2.5. Thus the fast components of $E^u(0,\xi)$ continued to rotate in the clockwise direction.

The conclusion matches the intuition that we have from Sturm-Liouville theory; the stable structures are the ones that oscillate less. The difference in this case is that the profiles of the stable and unstable waves look exactly the same. Instead, it is the the unstable bundles for these waves that do or do not oscillate. To connect this to phase space, the stable waves are the ones whose unstable manifolds don't contain any full twists as they go from $-\infty$ to $+\infty$.

\subsection{Discussion}
We will close by comparing the Maslov index calculation of this paper with that of \cite{CJ18}. The latter paper analyzed \emph{traveling} waves for the following FitzHugh-Nagumo system: \begin{equation}\label{FHN PDE}
\begin{aligned} u_t & =u_{zz}+f(u)-v\\
v_t & = v_{zz}+\eps(u-\gamma v).
\end{aligned}
\end{equation}
The traveling wave and eigenvalue equations are each four-dimensional for (\ref{FHN PDE}), hence the unstable bundle is a curve in $\Lambda(2)$. Nonetheless, beyond the difference in dimension, the two Maslov index calculations are quite similar. Most notably, the waves in both cases are homoclinic orbits consisting of two fast jumps between slow dynamics. Moreover, the Maslov index in each case comes down to a corner trajectory in the Lagrangian Grassmannian, as in \S \ref{sec corner}.

The key difference between the two models is the mechanics of the corner calculation. In \S \ref{sec corner}, the task was to find a trajectory that connected a fixed point with a regular point in $\Lambda(3)$. Such a trajectory is obviously unique. Conversely, the corner calculation for (\ref{FHN PDE}) amounts to finding a heteroclinic orbit connecting two fixed points, cf. \cite[\S 4.4]{CJ18}. Using the invariance of $\Lambda(2)$ for the eigenvalue equation, it is shown that two such orbits exist. Exactly one of the two possible paths crosses the singular cycle, and stability of the waves hinges on which path is taken.\footnote{This is true specifically for the \emph{second} corner, which is the transition from the first slow segment to the fast back. Although multiple heteroclinic connections exist at the other two corners as well, neither path taken would result in a conjugate point.} We showed that the wave is stable by keeping track of the orientation of the unstable manifold during the fast-slow transitions.

The technical reason for this difference is that the transversality condition needed to prove existence of the waves (via the Exchange Lemma) requires $\eps>0$ for (\ref{PDE model slow}), but not for (\ref{FHN PDE}). Indeed, the jump-off condition (\ref{jump-off_cond_dependet}) determines the set of points which land on $M_\eps^+$, and hence the orientation of $W^u(X_-)$ at $z_2$. The tangent space here is given by (\ref{def Z}), which is not a fixed point for (\ref{linearization_corner}). In the traveling wave case, varying the speed $c$ allows one to obtain transversality with $\eps=0$, so the reorientation of the unstable bundle at the corners is simply an exchange of eigenvectors.

Ignoring the technical reason, the two calculations give the impression that (in)stability is ``inevitable" in the standing wave case but ``fortuitous" in the traveling wave case. To be more precise, the fate of the standing wave (with respect to stability) is sealed along the first fast jump. For the traveling wave, on the other hand, the potential conjugate point lives in the second corner. More importantly, there is no obvious reason to select one heteroclinic orbit over the other for the reorientation leading into the fast back. One is tempted to see if the FitzHugh-Nagumo fast waves can be \emph{destabilized} by somehow altering the equations to select the other orbit instead.

The preceding discussion speaks to the motivation for using the Maslov index to study stability in the first place. By using detailed information about the wave and its phase space, we can obtain intrinsic reasons for stability or instability. Furthermore, this information can be used to distinguish different mechanisms for instability, and perhaps show us how to generate such instabilities.

\paragraph*{\textbf{Acknowledgments:}} C.J. acknowledges support from the US Office of Naval Research under grant number N00014-18-1-2204.

\bibliographystyle{amsplain}
\bibliography{refs}

\providecommand{\bysame}{\leavevmode\hbox to3em{\hrulefill}\thinspace}
\providecommand{\MR}{\relax\ifhmode\unskip\space\fi MR }
\providecommand{\MRhref}[2]{%
  \href{http://www.ams.org/mathscinet-getitem?mr=#1}{#2}
}
\providecommand{\href}[2]{#2}
\begin{thebibliography}{10}

\bibitem{AGJ}
J.~Alexander, R.A. Gardner, and C.K.R.T. Jones, \emph{A topological invariant
  arising in the stability analysis of travelling waves}, J. reine angew. Math
  \textbf{410} (1990), no.~167-212, 143.

\bibitem{Arnold67}
Vladimir~Igorevich Arnol'd, \emph{Characteristic class entering in quantization
  conditions}, Functional Analysis and its applications \textbf{1} (1967),
  no.~1, 1--13.

\bibitem{BCJLMS17}
Margaret Beck, Graham Cox, Christopher Jones, Yuri Latushkin, Kelly McQuighan,
  and Alim Sukhtayev, \emph{Instability of pulses in gradient
  reaction--diffusion systems: a symplectic approach}, Phil. Trans. R. Soc. A
  \textbf{376} (2018), no.~2117, 20170187.

\bibitem{BJ}
Amitabha Bose and Christopher~K.R.T. Jones, \emph{Stability of the in-phase
  travelling wave solution in a pair of coupled nerve fibers}, Indiana
  University Mathematics Journal \textbf{44} (1995), no.~1, 189--220.

\bibitem{Bru96}
Pavol Brunovsk{\`y}, \emph{Tracking invariant manifolds without differential
  forms}, Acta Math. Univ. Comenianae \textbf{65} (1996), no.~1, 23--32.

\bibitem{Bru99}
\bysame, \emph{${C}^r$-inclination theorems for singularly perturbed
  equations}, Journal of Differential Equations \textbf{155} (1999), no.~1,
  133--152.

\bibitem{CRS}
Paul Carter, Jens~DM Rademacher, and Bj{\"o}rn Sandstede, \emph{Pulse
  replication and accumulation of eigenvalues}, arXiv preprint arXiv:2005.11683
  (2020).

\bibitem{CDB09}
Fr{\'e}d{\'e}ric Chardard, Fr{\'e}d{\'e}ric Dias, and Thomas~J Bridges,
  \emph{Computing the {M}aslov index of solitary waves, {P}art 1: {H}amiltonian
  systems on a four-dimensional phase space}, Physica D: Nonlinear Phenomena
  \textbf{238} (2009), no.~18, 1841--1867.

\bibitem{CDB_multi}
\bysame, \emph{On the maslov index of multi-pulse homoclinic orbits},
  Proceedings of the Royal Society A: Mathematical, Physical and Engineering
  Sciences \textbf{465} (2009), no.~2109, 2897--2910.

\bibitem{CC12}
Chao-Nien Chen and YS~Choi, \emph{Standing pulse solutions to
  {F}itz{H}ugh--{N}agumo equations}, Archive for Rational Mechanics and
  Analysis \textbf{206} (2012), no.~3, 741--777.

\bibitem{CH}
Chao-Nien Chen and Xijun Hu, \emph{Maslov index for homoclinic orbits of
  {H}amiltonian systems}, Annales de l'Institut Henri Poincar{\'e}, vol.~24,
  Analyse non Lin{\'e}are, no.~4, 2007, pp.~589--603.

\bibitem{CH14}
\bysame, \emph{Stability analysis for standing pulse solutions to
  {F}itz{H}ugh--{N}agumo equations}, Calculus of Variations and Partial
  Differential Equations \textbf{49} (2014), no.~1-2, 827--845.

\bibitem{corn}
Paul Cornwell, \emph{Opening the {M}aslov {B}ox for traveling waves in
  skew-gradient systems: counting eigenvalues and proving (in) stability},
  Indiana University Mathematics Journal \textbf{68} (2019), no.~6, 1801--1832.

\bibitem{CJ18}
Paul Cornwell and Christopher~K.R.T. Jones, \emph{On the existence and
  stability of fast traveling waves in a doubly diffusive
  {F}itz{H}ugh--{N}agumo system}, SIAM Journal on Applied Dynamical Systems
  \textbf{17} (2018), no.~1, 754--787.

\bibitem{DvHK_ex}
Arjen Doelman, Peter Van~Heijster, and Tasso~J Kaper, \emph{Pulse dynamics in a
  three-component system: existence analysis}, Journal of Dynamics and
  Differential Equations \textbf{21} (2009), no.~1, 73--115.

\bibitem{Dum93}
Freddy Dumortier, \emph{Techniques in the theory of local bifurcations:
  {B}low-up, normal forms, nilpotent bifurcations, singular perturbations},
  Bifurcations and periodic orbits of vector fields, Springer, 1993,
  pp.~19--73.

\bibitem{Fen79}
Neil Fenichel, \emph{Geometric singular perturbation theory for ordinary
  differential equations}, Journal of Differential Equations \textbf{31}
  (1979), no.~1, 53--98.

\bibitem{hatcher}
Allen Hatcher, \emph{Algebraic topology}, Cambridge University Press, 2002.

\bibitem{Henry}
Dan Henry, \emph{Geometric theory of semilinear parabolic equations}, Lecture
  notes in mathematics, Springer-Verlag, Berlin, New York, 1981.

\bibitem{HLS16}
Peter Howard, Yuri Latushkin, and Alim Sukhtayev, \emph{The {M}aslov and
  {M}orse indices for {S}chr{\"o}dinger operators on $\mathbb{R}$}, Indiana
  University Mathematics Journal \textbf{67} (2018), no.~5, 1765--1815.

\bibitem{Jones_GSP}
Christopher Jones, \emph{Geometric singular perturbation theory}, Dynamical
  systems (1995), 44--118.

\bibitem{JLS17}
Christopher Jones, Yuri Latushkin, and Selim Sukhtaiev, \emph{Counting spectrum
  via the {M}aslov index for one dimensional $\theta$-periodic
  {S}chr{\"o}dinger operators}, Proceedings of the American Mathematical
  Society \textbf{145} (2017), no.~1, 363--377.

\bibitem{Jo88}
Christopher~K.R.T. Jones, \emph{Instability of standing waves for non-linear
  {S}chr{\"o}dinger-type equations}, Ergodic Theory and Dynamical Systems
  \textbf{8} (1988), no.~8*, 119--138.

\bibitem{JLM13}
Christopher~K.R.T. Jones, Yuri Latushkin, and Robert Marangell, \emph{The
  {M}orse and {M}aslov indices for matrix {H}ill's equations}, Spectral
  analysis, differential equations and mathematical physics: a festschrift in
  honor of Fritz Gesztesy’s 60th birthday \textbf{87} (2013), 205--233.

\bibitem{JK94}
C.K.R.T. Jones and N.~Kopell, \emph{Tracking invariant manifolds with
  differential forms in singularly perturbed systems}, Journal of Differential
  Equations \textbf{108} (1994), no.~1, 64--88.

\bibitem{Kuehn15}
Christian Kuehn, \emph{Multiple time scale dynamics}, vol. 191, Springer, 2015.

\bibitem{Milnor}
John~Willard Milnor, \emph{Morse theory}, Princeton university press, 1963.

\bibitem{RS93}
Joel Robbin and Dietmar Salamon, \emph{The {M}aslov index for paths}, Topology
  \textbf{32} (1993), no.~4, 827--844.

\bibitem{Sandstede02}
Bj{\"o}rn Sandstede, \emph{Stability of travelling waves}, Handbook of
  dynamical systems \textbf{2} (2002), 983--1055.

\bibitem{Schenk}
CP~Schenk, M~Or-Guil, M~Bode, and H-G Purwins, \emph{Interacting pulses in
  three-component reaction-diffusion systems on two-dimensional domains},
  Physical Review Letters \textbf{78} (1997), no.~19, 3781.

\bibitem{vHCNT}
Peter van Heijster, Chao-Nien Chen, Yasumasa Nishiura, and Takashi Teramoto,
  \emph{Localized patterns in a three-component {F}itz{H}ugh--{N}agumo model
  revisited via an action functional}, Journal of Dynamics and Differential
  Equations \textbf{30} (2018), no.~2, 521--555.

\bibitem{vHDK_stab}
Peter van Heijster, Arjen Doelman, and Tasso~J Kaper, \emph{Pulse dynamics in a
  three-component system: stability and bifurcations}, Physica D: Nonlinear
  Phenomena \textbf{237} (2008), no.~24, 3335--3368.

\end{thebibliography}

\end{document}